\documentclass[12pt]{article}
\usepackage{amsfonts,amsmath,amsxtra}
\usepackage{amssymb,amscd}
\usepackage[mathscr]{eucal}
\def\hybrid{\topmargin 0pt      \oddsidemargin 0pt
        \headheight 0pt \headsep 0pt
        \textwidth 16.5cm
        \textheight 23cm
        \marginparwidth 0.0in
        \parskip 5pt plus 1pt   \jot = 1.5ex}
\catcode`\@=11
\def\marginnote#1{}
\newcount\hour
\newcount\minute
\newtoks\amorpm
\hour=\time\divide\hour by60
\minute=\time{\multiply\hour by60 \global\advance\minute by-\hour}
\edef\standardtime{{\ifnum\hour<12 \global\amorpm={am}%
        \else\global\amorpm={pm}\advance\hour by-12 \fi
        \ifnum\hour=0 \hour=12 \fi
      \number\hour:\ifnum\minute<10 0\fi\number\minute\the\amorpm}}
\edef\militarytime{\number\hour:\ifnum\minute<10 0\fi\number\minute}

\def\draftlabel#1{{\@bsphack\if@filesw {\let\thepage\relax
   \xdef\@gtempa{\write\@auxout{\string
      \newlabel{#1}{{\@currentlabel}{\thepage}}}}}\@gtempa
   \if@nobreak \ifvmode\nobreak\fi\fi\fi\@esphack}
        \gdef\@eqnlabel{#1}}
\def\@eqnlabel{}
\def\@vacuum{}
\def\draftmarginnote#1{\marginpar{\raggedright\scriptsize\tt#1}}

\def\draft{\oddsidemargin -0.1truein
        \def\@oddfoot{\sl preliminary draft \hfil
        \rm\thepage\hfil\sl\today\quad\militarytime}
        \let\@evenfoot\@oddfoot \overfullrule 3pt
        \let\label=\draftlabel
        \let\marginnote=\draftmarginnote
\def\@eqnnum{{\rm (\theequation)}
\rlap{\kern\marginparsep\tt\@eqnlabel}%
\global\let\@eqnlabel\@vacuum}  }

\newcommand{\RR}{{\mathbb{R}}}
\newcommand{\CC}{{\mathbb{C}}}

\newcommand{\ZZ}{{\mathbb{Z}}}

\newfont{\Bbbb}{msbm7 scaled 1\@ptsize00}
\newcommand{\zs}{\raise-1pt\hbox{$\mbox{\Bbbb Z}$}}

\font\teneufm=cmmib10 scaled 1\@ptsize00
\font\seveneufm=cmmib7 scaled 1\@ptsize00
\font\fiveeufm=cmmib5  
\def\bfit#1{{\textfont1=\teneufm\scriptfont1=\seveneufm
\scriptscriptfont1=\fiveeufm
\mathchoice{
\hbox{$\mathsurround=0pt\displaystyle#1$}}
{\mathsurround=0pt\hbox{$\textstyle#1$}}
{\hbox{$\mathsurround=0pt\scriptstyle#1$}}
{\hbox{$\mathsurround=0pt\scriptscriptstyle#1$}}}}

\font\sevenmsa=msam6 
\newfam\msafam
\textfont\msafam=\sevenmsa
\def\hexnumber@#1{\ifnum#1<10 \number#1\else
\ifnum#1=10 A\else\ifnum#1=11 B\else\ifnum#1=12 C\else
\ifnum#1=13 D\else\ifnum#1=14 E\else\ifnum#1=15 F\fi\fi\fi\fi\fi\fi\fi}
\def\msa@{\hexnumber@\msafam}
\def\llcorner{\delimiter"4\msa@78\msa@78 }
\def\lrcorner{\delimiter"5\msa@79\msa@79 }
\mathchardef\blacktriangleright="3\msa@49
\mathchardef\blacktriangleleft="3\msa@4A
\font\tenmsb=msbm10 scaled 1\@ptsize00
\newfam\msbfam
\textfont\msbfam=\tenmsb
\scriptfont\msbfam=\tenmsb

\newdimen\linethick  \linethick=0.4pt
\newdimen\hboxitspace    \hboxitspace=5pt
\newdimen\vboxitspace    \vboxitspace=5pt

\def\fr#1{%
\be\new
\vcenter{
\hrule height\linethick
           \hbox{\vrule width\linethick
                 \kern\hboxitspace
                 \vbox{\kern\vboxitspace
                       \hbox{$\begin{array}{c}\displaystyle#1
          \end{array}$}%
                       \kern\vboxitspace}%
                 \kern\hboxitspace
                 \vrule width\linethick}%
           \hrule height\linethick}%
\ee}

\newdimen\Squaresize \Squaresize=14pt
\newdimen\Thickness \Thickness=0.5pt

\def\Square#1{\hbox{\vrule width \Thickness
   \vbox to \Squaresize{\hrule height \Thickness\vss
      \hbox to \Squaresize{\hss#1\hss}
   \vss\hrule height\Thickness}
\unskip\vrule width \Thickness}
\kern-\Thickness}

\def\Vsquare#1{\vbox{\Square{$#1$}}\kern-\Thickness}

\def\numberbysection{\@addtoreset{equation}{section}
        \def\theequation{\thesection.\arabic{equation}}}
\numberbysection

\renewcommand{\theequation}{\thesection.\arabic{equation}}
\def\titlepage{\@restonecolfalse\if@twocolumn\@restonecoltrue\onecolumn
     \else \newpage \fi \thispagestyle{empty}\c@page\z@
        \def\thefootnote{\fnsymbol{footnote}} }

\def\endtitlepage{\if@restonecol\twocolumn \else  \fi
        \def\thefootnote{\arabic{footnote}}
        \setcounter{footnote}{0}}  
\relax

\hybrid
\makeatletter
\newdimen\normalarrayskip            
\newdimen\minarrayskip               
\normalarrayskip\baselineskip
\minarrayskip\jot
\newif\ifold             \oldtrue            \def\new{\oldfalse}
\def\arraymode{\ifold\relax\else\displaystyle\fi}
\def\eqnumphantom{\phantom{(\theequation)}} 
\def\@arrayskip{\ifold\baselineskip\z@\lineskip\z@
     \else
     \baselineskip\minarrayskip\lineskip1\baselineskip\fi}

\def\@arrayclassz{\ifcase \@lastchclass \@acolampacol \or
\@ampacol \or \or \or \@addamp \or
   \@acolampacol \or \@firstampfalse \@acol \fi
\edef\@preamble{\@preamble
  \ifcase \@chnum
     \hfil$\relax\arraymode\@sharp$\hfil
     \or $\relax\arraymode\@sharp$\hfil
     \or \hfil$\relax\arraymode\@sharp$\fi}}

\def\@array[#1]#2{\setbox\@arstrutbox=\hbox{\vrule
     height\arraystretch \ht\strutbox
     depth\arraystretch \dp\strutbox
width\z@}\@mkpream{#2}\edef\@preamble{\halign \noexpand\@halignto
\bgroup \tabskip\z@ \@arstrut \@preamble \tabskip\z@ \cr}%
\let\@startpbox\@@startpbox \let\@endpbox\@@endpbox
  \if #1t\vtop \else \if#1b\vbox \else \vcenter \fi\fi
  \bgroup \let\par\relax
  \let\@sharp##\let\protect\relax
  \@arrayskip\@preamble}
\def\eqnarray{\stepcounter{equation}%
              \let\@currentlabel=\theequation
              \global\@eqnswtrue
              \global\@eqcnt\z@
              \tabskip\@centering              
              \let\\=\@eqncr
              $$%
            \halign to \displaywidth  \bgroup
             \eqnumphantom \@eqnsel
      \hskip\@centering                               
    $\displaystyle  \tabskip\z@ {##}$%
    &\global\@eqcnt\@ne \hskip 2\arraycolsep
         $ \displaystyle  \arraymode{##}$\hfil
    &\global\@eqcnt\tw@ \hskip 2\arraycolsep
         $\displaystyle\tabskip\z@{##}$\hfil
         \tabskip\@centering
    &{##}\tabskip\z@\cr}
\makeatother
\newtheorem{te}{Theorem}[section]

\newtheorem{prop}{Proposition}[section]           
\newtheorem{cor}{Corollary}[section]
\newtheorem{lem}{Lemma}[section]

\newcommand{\beq}[1]{\begin{equation}\label{#1}}
\newcommand\eeq{\end{equation}}
\newcommand\bqa{\begin{eqnarray}}
\newcommand\eqa{\end{eqnarray}}
\def\be{\begin{eqnarray}\new\begin{array}{cc}}
\def\ee{\end{array}\end{eqnarray}}

\def\beq{\begin{equation}}
\def\eeq{\end{equation}}
\def\bse{\begin{subequations}}                
\def\ese{\end{subequations}}
\def\bp{\begin{pmatrix}}
\def\ep{\end{pmatrix}}

\def\h{\hbar}
\def\i{\imath}

\def\square{\hfill{\vrule height6pt width6pt            
depth1pt} \break \vspace{.01cm}}                        

\def\stack#1#2{\raise0.7pt\hbox{$\mathrel{\mathop{#2}\limits^{#1}}$}}
\def\tr{\triangleright}
\def\tl{\triangleleft}
\def\sem{\mathsurround=0pt \raise1pt
\hbox{$\scriptscriptstyle>\!\!$}\:\!\!\tl}
\def\mes{\mathsurround=0pt \tr\!\:\!\raise0.8pt
\hbox{$\scriptscriptstyle\!\!<$}\,}
\def\]{\mathsurround=0pt ]\raise-2pt\hbox{$_\ast$}}
\def\balpha{{\bfit\alpha}}

\def\bmu{{\bfit\mu}}

\def\blambda{{\bfit\lambda}}
\def\brho{{\bfit\rho}}

\def\la{\lambda}

\def\<{\langle}
\def\>{\rangle}

\def\frak{\mathfrak}

\def\we{\raise-1pt\hbox{$\,\stackrel{\wedge}{,}\,$}}

\def\pr {\partial}

0
\begin{document}

\footnotesize
\normalsize

\newpage

\thispagestyle{empty}

\begin{center}

\phantom.
\bigskip
{
\hfill{\normalsize ITEP-TH-34/05}\\
\hfill{\normalsize HMI-05-11}\\
[10mm]\Large\bf
On a Gauss-Givental Representation of Quantum Toda

Chain Wave Function}
\vspace{0.5cm}

\bigskip\bigskip
{\large A. Gerasimov}
\\ \bigskip
{\it Institute for Theoretical \& Experimental Physics, 117259,
Moscow, Russia}\\ {\it
 Department of Pure and Applied Mathematics, Trinity
College, Dublin 2, Ireland } \\ {\it Hamilton
Mathematics Institute, TCD, Dublin 2, Ireland}\\
{\it Max-Planck-Institut fur Mathematik,
Vivatsgasse 7, D-53111
Bonn, Germany},\\
\bigskip
{\large S. Kharchev\footnote{E-mail: kharchev@itep.ru}},\\
\bigskip
{\it Institute for Theoretical \& Experimental Physics, 117259,
Moscow,
Russia}\\
\bigskip
{\large D. Lebedev\footnote{E-mail: lebedev@mpim-bonn.mpg.de}}
\\ \bigskip
{\it Institute for Theoretical \& Experimental Physics, 117259,
Moscow, Russia}  {\it Max-Planck-Institut fur Mathematik,
Vivatsgasse 7, D-53111
Bonn, Germany},\\
\bigskip
{\large S. Oblezin} \footnote{E-mail: Sergey.Oblezin@itep.ru}\\
\bigskip {\it Institute for Theoretical \& Experimental Physics,
117259, Moscow,
Russia}\\

\end{center}

\vspace{0.5cm}

\begin{abstract}
\noindent

We propose group theory interpretation of the integral
representation of the quantum open Toda chain wave function due to
Givental. In particular we construct  the  representation of
$U(\mathfrak{gl}(N))$ in terms of first order differential operators
in Givental variables. The construction of this  representation turns out
to be
closely  connected  with the integral representation based on the
factorized Gauss decomposition. We also reveal the recursive structure
of the Givental representation and provide the connection
with the Baxter $Q$-operator formalism.  Finally the generalization of the
integral representation to  the infinite and periodic quantum Toda wave functions
is discussed.
\end{abstract}

\vspace{1cm}

\clearpage \newpage


\normalsize
\section{Introduction}

In 1996 A. Givental proposed a new  integral formula for  the common
eigenfunction of open Toda chain  Hamiltonian operators  \cite{Gi} (see also
\cite{JK}). The first non trivial Hamiltonian operator of the $N$ particle
open Toda chain is given by
\be
H=-\frac{\hbar^2}{2}\sum\limits_{i=1}^{N}\frac{\partial^2}{\partial x_i^2}+
\sum\limits_{i=1}^{N-1}e^{x_{i}-x_{i+1}}\,.
\ee
The new integral formula for the common eigenfunction is given by the
following multiple integral
\cite{Gi}, \cite{JK}:
\be\label{Givp}
\Psi_{\la_1,\ldots,\la_N}(T_{N,1},\ldots,T_{N,N})=\int_{\Gamma}
e^{\frac{1}{\hbar}\mathcal{F}_N(T)}\prod_{k=1}^{N-1}
\prod_{i=1}^k dT_{k,i},
\ee
where $T_{N,i}:=x_i$, the function $\mathcal{F}_N(T)$ is given by
\be\label{pot}
\hspace{-0.5cm}
\mathcal{F}_N(T)={\imath}\sum_{k=1}^N\lambda_k\Big(\sum_{i=1}^kT_{k,i}-
\sum_{i=1}^{k-1}T_{k-1,i}\Big)-
\sum_{k=1}^{N-1}\sum_{i=1}^k\Big(e^{T_{k+1,i}-T_{k,i}}
+e^{T_{k,i}-T_{k+1,i+1}}\Big),
\ee
and the cycle $\Gamma$ is a middle dimensional submanifold in the $N(N-1)/2$-
dimensional complex torus with coordinates
$\{\exp\,{T_{k,i}},\, i=1,\ldots,k;\,k=1,\ldots, N-1\}$ such that the integral
converges. In particular one can choose $\Gamma=\RR^{N(N-1)/2}$. This integral
representation was motivated by the constructions of the Quantum Cohomology
of the flag manifolds and its description in mirror dual Landau-Ginzburg
model.

Basically  two different regular approaches to the problem of the explicit
construction of the integral representations of the open Toda chain wave
function are known. The first goes back to B. Kostant (1978) and reduces the
eigenvalue problem to the construction of a particular matrix element of an
irreducible representation of the corresponding Lie algebras \cite{Ko1},
\cite{STS}. Several integral representations of this matrix element are known
\cite{J}, \cite{SCH}, \cite{St}, \cite{GKMMO}. This approach may be
generalized to the periodic Toda chain using the representation theory of the
affine Lie groups but has not yet lead to the explicit integral formulas for
the wave functions in the periodic case.

Another approach is based  on  the Quantum Inverse Scattering Method (QISM)
\cite{Fa}, \cite{KS}, applied to both open and periodic quantum Toda chains
in \cite{Gu}, \cite{Ga}, \cite{Sk1}. The explicit integral representations
for the eigenfunctions in this framework were constructed in \cite{KL1},
\cite{KL2}, \cite{KL3}. The explanation  of this representation for the open
Toda chain in terms of the representation theory of Lie groups was given in
\cite{GKL1} (see also \cite{GKL2}) and was based on the generalization of the
Gelfand-Zetlin construction \cite{GelZ}, \cite{GelG} to the case of the
infinite-dimensional representations of $U(\mathfrak{gl}(N))$ proposed in
\cite{GKL1}.

Thus the representation theory provides a unifying framework for the
constructions of the integral representations of the wave functions of open
Toda chain. This naturally raises the question of the interpretation of the
integral formula (\ref{Givp}), (\ref{pot}) in terms of the representation
theory.

In this note we demonstrate that the explicit integral representation
(\ref{Givp}), (\ref{pot}) for $N$ particle open Toda chain naturally arises in
representation theory approach. It is based on a particular parameterization
of the upper/lower triangular parts of the Gauss decomposition of the group
element of $GL(N,\mathbb{R})$. This parameterization is close to the
factorization into the product of the elementary Jacobi matrices (see e.g.
\cite{BFZ}) but is slightly different. An interesting property of this
parameterization is that it provides a simple construction of the principal
series representation of the  universal enveloping algebra
$U(\mathfrak{gl}(N))$ acting in the space of functions on the {\it totally
positive} unipotent upper-triangular matrices. The representation is given by
the following explicit expressions for the images
$E^{(N)}_{i,j}$, $i,j=1,\ldots, N$ of the generators of $\frak{gl}(N)$ as the
first order differential operators in variables $T_{k,i}$, $i=1,\ldots,k;\,
k=1,\ldots, N-1$
\be\label{rep}
E^{(N)}_{i,i}=\mu^{(N)}_i+ \sum_{k=1}^{i-1}\frac{\partial}{\partial
T_{N+k-i,k}}-\sum_{k=i}^{N-1}\frac{\partial}{\partial T_{k,i}}\,,\\
E^{(N)}_{i,i+1}=\sum_{n=1}^i
\left(\,\sum_{k=n}^ie^{T_{N+k-i,k}-T_{N+k-i-1,k}}\right)\!\!
\left(\frac{\partial}{\partial T_{N+n-i-1,n}}-
\frac{\partial}{\partial T_{N+n-i-1,n-1}}\right),\\
E^{(N)}_{i+1,i}=\sum_{k=i}^{N-1}e^{T_{k,i}-T_{k+1,i+1}}
\left(\mu^{(N)}_i-\mu^{(N)}_{i+1}+\sum_{s=i}^k\left(\frac{\partial}{\partial
T_{s,i+1}}-\frac{\partial}{\partial T_{s,i}}\right)\right),
\ee
where $\mu^{(N)}_i=-\i \h^{-1}\la_i\,-\,\rho_i^{(N)}$, $\la_i\in\RR^N$,
$\rho_i^{(N)}=\frac{1}{2}(N-2i+1),\,i=1,\ldots N$ and we assume that
$T_{N,i}=0$ and $\sum_{i}^j=0$  for $i>j$. The functions of the
variables $T_{k,i}$ are identified with the functions
on the subspace of the totally positive unipotent upper-triangular matrices.
In the following we will use  the term Gauss-Givental representation for
this representation and the resulted representation of the open Toda wave
function.

Let us stress that the Gauss-Givental representation is closely related with
the recursive structure of the open Toda chains \cite{Sk2}. In particular
there exists the explicit integral recursion operator connecting the wave
functions of the $N-1$ particle and the $N$ particle Toda chains. The
iterative application of this operator provides a simple independent
derivation of the integral representation \cite{Gi}, \cite{JK}. Note that the
recursion operator is  closely connected with the Baxter $Q$-operator for
periodic Toda chain \cite{Ga}, \cite{PG}. Thus from this point of view,
the Gauss-Givental representation is a direct consequence of the generalized
Baxter $Q$-operator formalism. This provides another relation between
representation theory and QISM approaches to the solution of open Toda chain.
Note that the similar iterative procedure also based on Baxter $Q$-operator
was used in \cite{GKL1} to get another integral representation in
Gelfand-Zetlin parameterization. Finally let us remark that in \cite{St} the
iterative construction of the Whittaker function was proposed in the framework
of the Iwasawa decomposition of the group element. It is easy to see that
after simple manipulations this iterative procedure may be reduced to the
one discussed in this paper.

The plan of the paper is as follows. In Section 2 we construct a
representation of $U(\mathfrak{gl}(N))$ in terms of differential operators
using a particular parameterization of the totally positive upper-triangular
matrices. In Section 3 we derive the integral representation for the Toda
chain wave function using this parameterization and show that integral
representation obtained in this way coincides with the representation proposed
by Givental. In Section 4 we rederive the integral representation using the
iterative procedure and discuss the connection with Baxter $Q$-operator
formalism. In Section 5, using the results of the previous section, we
consider a generalization of the integral representations of wave function to
infinite and periodic Toda chains. We conclude in Section 6 with the
discussions of further directions of research. In the Appendices A and B the
proofs of various statements  from the main part of the text are given.

{\em Acknowledgments}:  The authors are grateful to M. Finkelberg and
A. Rosly for useful discussions. The research was partly supported by grants
CRDF RM1-2545; INTAS 03-513350; grant NSh 1999.2003.2 for support of
scientific schools, and also by grants RFBR-03-02-17554 (A. Gerasimov,
D. Lebedev, S. Oblezin), and RFBR-03-02-17373 (S. Kharchev).
The research of A. Gerasimov was also partly supported by SFI Basic
Research Grant.  S. Oblezin would like to thank l'Institute des Hautes
\'Etudes Scientifiques and the Max-Planck-Institut f\"{u}r Mathematik for
support and warm hospitality. The research of S. Oblezin was also
partly supported as a Independent University M\"obius Prize Fellow.

\section{Representation of $U(\mathfrak{gl}(N))$}

In this section we construct a particular realization of the
principal series representation of $U(\mathfrak{g}(N))$ acting in the
space of  functions  on  the subset of the totally positive
elements of the nilpotent subgroup of $G=GL(N,\mathbb{R})$.

Let $N_+$  be a unipotent radical of the  Borel subgroup $B_+$ of the group
$G$. We fix the choice of this subgroups by considering the elements of $B_+$
as upper triangular matrices and the elements of $N_+$ as upper triangular
matrices with unit diagonal. We define the subset $N^{(+)}_+\subset N_+$ of
the totally  positive unipotent upper-triangular matrices. The matrix
$x\in N_+$ is called totally positive if all its not identically zero minors
are positive real numbers (for details see e.g. \cite{BFZ}). Consider the
subspace of functions $M_{\bmu}\subset M=Fun(G)$ equivariant with
respect to the left action of the opposite Borel subgroup of the
lower-triangular matrices $B_-\subset G$: $f(bg)=\chi _{\bmu}(b)f(g)$,
where $b\in B_-$, $\bmu=(\mu^{(N)}_1,\ldots,\mu^{(N)}_N)$ and
$\chi _{\bmu}(b)=\prod\limits_{i=1}^{N} |b_{ii}|^{\mu^{(N)}_i}$ is a
character of $B_-$.

The right action of the group on  $M_{\bmu}$ defines a realization
of the principal series representation of $GL(N,\mathbb{R})$  for
$\mu_i^{(N)}=-\i\h^{-1}\la_i-\rho_i^{(N)}$, where $\la_i\in\RR$ and
$\rho_i^{(N)}=\frac{1}{2}(N-2i+1)$, $i=1,\ldots N$ are the
components of the vector $\brho:=\frac{1}{2}\sum_{\balpha>0}\balpha$
in the standard basis of $\RR^N$. Let $M^{(+)}_{\bmu}$ by the set of
functions on the totally positive unipotent upper-triangular
matrices. We will consider the representation of the universal
enveloping algebra $U(\mathfrak{gl}(N))$ in $M^{(+)}_{\bmu}$ induced
by the above construction. Let $e_{i,j}$ stand for the elementary
$N\times N$ matrix with the unit at the $(i,j)$ place and zeros
otherwise. Consider the set of the diagonal elements parameterized
as follows
$U_k=\sum\limits_{i=1}^{k}e^{T_{k,i}}e_{i,i}+\sum\limits_{i=k+1}^{N}
e_{i,i}$. Define the following set of the upper-triangular matrices
\be\label{deformed} \tilde{U}_k=\sum\limits_{i=1}^{k}
e^{T_{k,i}}e_{i,i} + \sum\limits_{i=k+1}^{N} e_{i,i} +
\sum\limits_{i=1}^{k-1} e^{T_{k-1,i}}e_{i,i+1}. \ee Then any totally
positive element  $x\in N^{(+)}_+$ can be represented in the
following form: \be\label{givpar}
x=\tilde{U}_{2}U_{2}^{-1}\tilde{U}_{3}U_{3}^{-1}\cdots
\tilde{U}_{N-1}U_{N-1}^{-1}\tilde{U}_{N},
\ee where we assume that $T_{N,i}=0$. This can be easily verified
using the connection with the parametrization of the upper
triangular matrices in terms of the product of the elementary Jacobi
matrices  (see e.g. \cite{BFZ}),
 namely, an arbitrary element $x\in N^{(+)}_+$ can be represented  as
follows:
\bqa\label{factor}
x= \prod_{k=1}^{N-1}\left(1+\sum_{i=1}^k
y_{k,i}e_{i,i+1}\right),
\eqa
where the variables $y_{k,i}$ are positive numbers. The comparison of
(\ref{givpar}), (\ref{factor}) leads to the following relations:
\bqa\label{ch1}
y_{k,i}=e^{T_{k,i}-T_{k+1,i+1}}\,.
\eqa
Now let us use this parameterization to construct a representation of
$U(\mathfrak{gl}(N))$ in $M^{(+)}_{\mu}$.
The images, $E_{i,j}^{(N)}$, of the generators of $U(\mathfrak{gl}(N))$ in the representation induced
from the one dimensional representation of the Borel subalgebra are
defined as follows. Let $E^{(N)}_{i,j}$ be a first order differential
operator in $T_{k,i}$ such that for  any  $f\in M_{\bmu}$
one has
 \bqa\label{def}
E^{(N)}_{i,j}f(x)=\frac{d}{d \epsilon}f(x(1+\epsilon
e_{i,j}))\!\!\mid_{\,\epsilon=0}.
\eqa
\begin{prop}\label{baserep}
The following differential operators define a 
representation $\pi_{\bmu}$ of $\,\frak{gl}(N)$ in $M_{\bmu}$:
\be
E^{(N)}_{i,i}=\mu^{(N)}_i+ \sum_{k=1}^{i-1}\frac{\partial}{\partial
T_{N+k-i,k}}-\sum_{k=i}^{N-1}\frac{\partial}{\partial
T_{k,i}}\,,\,\\
E^{(N)}_{i,i+1}=\sum_{n=1}^i\left(\,
\sum_{k=n}^ie^{T_{N+k-i,k}-T_{N+k-i-1,k}}\right)
\!\!\left(\frac{\partial}{\partial T_{N+n-i-1,n}}-
\frac{\partial}{\partial T_{N+n-i-1,n-1}}\right),\\
E^{(N)}_{i+1,i}=\sum_{k=i}^{N-1}e^{T_{k,i}-T_{k+1,i+1}}
\left(\mu^{(N)}_i-\mu^{(N)}_{i+1}+\sum_{s=i}^k\left(\frac{\partial}{\partial
T_{s,i+1}}-\frac{\partial}{\partial T_{s,i}}\right)\right),
\ee
where we assume that $T_{N,i}=0$ and omit the derivatives over $T_{i,j},\,
i<j$.
\end{prop}
The proof is given in the Appendix A.

\section{Integral representation of the wave function}

\subsection{Matrix elements as open Toda wave functions}

In this section we construct an  integral representation of the wave
function of the open Toda chain  corresponding to the representation  $\pi_{\bmu}$
of $U(\mathfrak{gl}(N))$ introduced in the previous section. The
construction can be done using either Gauss or  Iwasawa decompositions
of the group. In the paper we will use the Gauss decomposition and the
corresponding Whittaker model of representation.

We first recall some facts from \cite{Ko2}. Let $\mathfrak{n}_+$ and
$\mathfrak{n}_-$ be two nilpotent  subalgebras of $\frak{gl}(N)$ generated,
respectively,  by positive and negative root generators.
The homomorphisms (characters) $\chi_+\!:n_+\!\rightarrow\CC$,
$\chi_-\!:n_-\!\!\rightarrow\CC$ are uniquely determined by their values on
the simple root generators, and are called non-singular if the complex numbers
$\chi_+(E^{(N)}_{i,i+1}):=\xi^{(i)}_R$ and
$\chi_-(E^{(N)}_{i+1,i}):=\xi^{(N-i)}_L$ are non-zero for all
$i=1,\ldots,N-1$.

Let $V$ be any $U(\frak{gl}(N))$-module.
A vector $\psi_R\in V$ is called a Whittaker vector with respect to the
character $\chi_+$ if
\be\label{ewt}
\hspace{2cm}
E^{(N)}_{i,i+1}\psi^{(N)}_R=\xi^{(i)}_R \psi^{(N)}_R\,,
\hspace{1cm}(i=1,\ldots,N-1),
\ee
and an element $\psi_L\in V'$ is called a Whittaker vector
with respect to the character $\chi_-$  if
\be\label{fwt}
\hspace{2cm} E^{(N)}_{i+1,i} \psi^{(N)}_L=\xi^{(N-i)}_L \psi^{(N)}_L\,,
\hspace{1cm} (i=1,\ldots,N-1).
\ee

A Whittaker vector is called cyclic in  $V$ if $U(\frak{gl}(N))\psi=V$,
and a $U(\frak{gl}(N))$-module is a Whittaker module if it contains a cyclic
Whittaker vector. The $U(\frak{gl}(N))$-modules $V$ and $V'$ are called dual
if there exists a non-degenerate pairing $\<.\,,.\>: V'\times V\to\CC$
such that $\<v',Xv\>=-\<X v',v\>$ for all $v\in V$, $v'\in V'$ and
$X\in\frak{gl}(N)$. Let us assume that the action of the Cartan subalgebra is
integrated to the action of the Cartan torus. An  eigenfunction of the
open Toda chain can  be written in terms of the pairing as follows
\bqa\label{pairing}
\Psi_{\blambda}^{(N)}(T_{N,1},\ldots,T_{N,N})=e^{-\sum
T_{N,i}\rho_i^{(N)}}\<\psi^{(N)}_L\,, \pi_{\bmu}(e^{-\sum T_{N,i}
E^{(N)}_{ii}}) \,\psi^{(N)}_R\>\,,
\eqa
where $\mu_k^{(N)}=-\i\h^{-1}\la_k-\rho_k^{(N)}$ and
$\rho_k^{(N)}=\frac{1}{2}(N-2k+1),\,k=1,\ldots N$
are the components of the vector $\brho:=\frac{1}{2}\sum\limits_{\balpha>0}\balpha$
 in the standard basis of $\RR^N$. Let $H_k$ be  the radial projections of
 generators $c_k$ of the center of the universal enveloping algebra
${\cal U}(\frak{gl}(N))$, defined by
\bqa \label{hamdef}
H_k \Psi^{(N)}_{\blambda}(T_{N,1},\ldots,T_{N,N})=e^{-\sum
T_{N,i}\rho_i^{(N)}}\<\psi^{(N)}_L\,,\pi_{\bmu}(e^{-\sum T_{N,i}
E^{(N)}_{ii}})\,c_k\,\psi^{(N)}_R\>.
\eqa
One can show that $H_k$ provide the complete set of commuting Hamiltonians of
the open Toda chain model.

\subsection{Matrix element in Gauss-Givental parameterization}

We are going to write down the  matrix element (\ref{pairing}) explicitly
in the Gauss-Givental representation introduced above.  Let us start with
construction of the Whittaker vectors $\psi_L^{(N)}$ and $\psi_R^{(N)}$.
\begin{prop}\label{Wvec}
The functions
\bqa\label{wR}
\psi^{(N)}_R=
\exp\left\{\sum_{i=1}^{N-1}\xi^{(i)}_R\sum_{k=i}^{N-1}
e^{T_{k,i}-T_{k+1,i+1}}\right\},\eqa and \bqa\label{wL}
\psi^{(N)}_L= \exp\left\{\sum_{k=1}^{N-1}\sum_{i=1}^k
(\mu_k-\mu_{k+1})T_{k,i}+
\sum_{i=1}^{N-1}\xi^{(i)}_{L}\sum_{k=1}^{N-i}e^{T_{k+i,k}-T_{k+i-1,k}}
\right\}, \eqa are solutions of the  linear differential equations
(\ref{ewt}) and (\ref{fwt}) (we assume that $T_{N,i}=0$).
\end{prop}
{\it Proof.}  The proof is based on  the recursion properties of the
Whittaker vectors $\psi^{(N)}_{L,R}$ and is  given in the Appendix
B.

  Evidently, the Whittaker modules $V:=U(\frak{gl}(N))\psi_R^{(N)}$ and
$V':=U(\frak{gl}(N))\psi_L^{(N)}$ are spanned by the elements
$\prod_{k=1}^{N-1}\prod_{i=1}^ke^{n_{k,i}T_{k,i}}\,\psi_{R,L}^{(N)}$, where
$n_{k,i}\in\ZZ$.

Now we construct the pairing between $V$ and $V'$.
 Define a measure of integration over $N_+^{(+)}$  as
$\omega_N=\wedge_{k=1}^{N-1}\wedge_{i=1}^ke^{T_{k,i}}d T_{k,i}$,
and for any $\phi\in V',\,\psi\in V$ introduce the following pairing:
\be\label{gpair}
\<\phi,\psi\>=\int\limits_\Gamma\omega_N\,\overline\phi(T)\psi(T)\,,
\ee
where the integration contour $\Gamma$ is chosen in such a way that the
integral (\ref{gpair}) is convergent for $\phi=\psi_L^{(N)}$ and
 $\psi=\psi_R^{(N)}$. Then the convergence for arbitrary elements
 $\phi \in V'$ and $\psi \in V$ follows. 
\begin{lem}
The pairing  (\ref{gpair}) is non-degenerate and satisfies the following
condition
\be \label{pcon}
\<\phi,X\psi\>=-\<X\phi,\psi\>\,,
\ee
for any $X\in \mathfrak{gl}(N)$.
\end{lem}

Now we are ready to find  the integral representation for the pairing 
(\ref{pairing}). To get the explicit expression for the integrand one should
make  a choice.  For example, the group element
$e^{H^{(N)}}=\exp(-\sum_{i=1}^NT_{N,i}E_{i,i}^{(N)})$
 can act on the right vector $\psi_R^{(N)}$ or on the left vector $\psi^{(N)}_L$.
We use the most symmetric choice.  Let us represent the Cartan group element
in the following way
\be\label{Cartan} e^{H^{(N)}}\,=\,e^{H^{(N)}_L}e^{H^{(N)}_R}=\\
\hspace{-0.5cm}
=\exp\left\{\sum_{i=1}^{N-1}
T_{N,i}\sum_{k=1}^{N-1}\frac{\pr}{\pr T_{k,i}}\right\}
 \exp\left\{-\sum_{i=1}^N \mu_i^{(N)}T_{N,i} -\sum_{i=2}^N
T_{N,i}\sum_{k=1}^{i-1}\frac{\pr}{\pr T_{N-i+k,k}}\right\}.
\ee
Then we suppose that $H^{(N)}_L$ acts on the left vector and $H^{(N)}_R$ acts
on the right vector in (\ref{pairing}). Taking into account the results of
the last  lemma one obtains the following theorem.
\begin{te}
The eigenfunction of the open Toda chain  defined by (\ref{pairing}),
(\ref{gpair}) has the following integral representation
\be\label{giv}
\Psi_{\blambda}^{(N)}(T_{N,1},\ldots,T_{N,N})=
\int\limits_\Gamma
\exp\left\{\frac{\i}{\hbar}\sum_{k=1}^N
\lambda_k\Big(\sum_{i=1}^kT_{k,i}-\sum_{i=1}^{k-1}T_{k-1,i}\Big)\right\}
\times\\
\exp\left\{\sum_{i=1}^{N-1}\bar\xi_{L}^{(i)}\sum_{k=1}^{N-i}
e^{T_{k+i,k}-T_{k+i-1,k}}+
\sum_{i=1}^{N-1}\xi_{R}^{(i)}\sum_{k=i}^{N-1}
e^{T_{k,i}-T_{k+1,i+1}}\right\}\prod_{k=1}^{N-1}\prod_{i=1}^kdT_{k,i}\,.
\ee
\end{te}
The integral representation (\ref{giv}) coincides with (\ref{Givp}),
(\ref{pot}), if we set $\xi^{(i)}_R=\xi^{(i)}_L=-\h^{-1}$.  In this case
one can choose the contour to be $\Gamma=\RR^{N(N-1)/2}$.

\subsection{Connection with another parametrization}

It is useful to compare the  integral representation obtained above with
 one  in \cite{GKMMO}, which is  based on a more standard
parameterization of the group element. Let us recall this construction.
We parameterize the unipotent upper triangular matrices $\|x_{i,j}\|$ by
its elements $x_{i,j}$,  $1\leq i< j\leq N$. In the corresponding  realization
of the representation of $U(\mathfrak{gl}(N))$ the action of the generators
of the algebra is given by the differential operators in $x_{i,j}$. The
measure of  integration entering the definition of  the  pairing between two
vectors in this realization is given by $$\omega'_N=
\wedge_{i=1}^{N-1}\wedge_{j=i+1}^Ndx_{i,j}.$$

Let us find the right Whittaker vector in this realization.
Integrating the action of the algebra $\mathfrak{n}_+$ to the action of the
group $N_+$ the equation for the right Whittaker vector can  be written in
the form
$$\pi _{\bmu}(z)\psi^{(N)}_{R}(x)=e^{\sum_{i=1}^{N-1}\xi_R^{(i)} z_{i,i+1}}
\psi^{(N)}_{R}(x),$$
for any $z\in N_+$. Thus we should find the one-dimensional representation of
the group of upper-tri\-angu\-lar matrices. Note that
the additive character of this group is easily constructed using the
fact that  the main diagonal behave additively under the multiplication.
Exponentiating  this character we get a one-dimensional representation:
\be\label{right}
\psi^{(N)}_{R}(x)=\prod_{i=1}^{N-1}e^{\xi_R^{(i)}x_{i,i+1}}.
\ee
Similarly, the left vector should satisfies the condition
$$
\pi _{\bmu}(z^t)\psi^{(N)}_{L}(x)= e^{\sum_{i=1}^{N-1}\xi_L^{(i)} z_{i,i+1}}
\psi^{(N)}_{L}(x),
$$
for any $z\in N_+$.
To construct such vector  we use an  inner automorphism
which maps the unipotent upper-triangular matrices to the unipotent
lower-triangular matrices. Consider the matrix
$(w_0)_{i,j}=\delta_{i+j,N+1}$ acting as  $z\to w_0^{-1}zw_0$. 
We have  $w_0^{-1}N_+w_0=N_-$ where $N_-$ is opposed to $N_+$.
 It is easy to check that we can define $\psi_L^{(N)}(x)=\psi_R^{(N)}(xw_0^{-1})$.
An  explicit calculation gives the following result
\be\label{left1}
\psi^{(N)}_{L}(x)=\prod_{i=1}^{N-1}\Big\{\Delta _i (x w_0^{-1})^
{\mu_i -\mu_{i+1}}e^{\xi_L^{(i)}
\frac{\Delta _{i,i+1}(x w_0^{-1})}{\Delta _{i}(xw_0^{-1})}}\Big\}\,,
\ee
where  $\Delta_i (M)$ denotes the principle $i\times i$ minor of the  matrix
$M$ and $\Delta_{i,i+1}(M)$ denotes the determinant obtained from
$\Delta_i (M)$ by interchanging the $i$-th and $(i+1)$-th columns in $M$.

Consider the decomposition (\ref{factor}) of a totally  positive
unipotent upper-triangular matrix $x$ where we   assume that $y_{i,j}$ are defined through
$T_{i,j}$ by the relation 
\be\label{ch1new} y_{k,i}=e^{T_{k,i}-T_{k+1,i+1}-T_{N,i}+T_{N,i+1}}.\ee

It is worth comparing the parameterizations (\ref{ch1}) and (\ref{ch1new}). 
 The  difference between the two reflects the fact that 
in  \cite{GKMMO} the group element
entering (\ref{pairing}) acts on the {\em right} vector in the integral
representation. This  difference may be  compensated by the shift of the
 integration variables $T_{k,i}\rightarrow T_{k,i}-T_{N,i}$ thus
 reconciling (\ref{ch1}) and (\ref{ch1new}).  

Following \cite{BFZ}, one can
express the minors through the variables $y_{k,i}$ as well as $T_{k,i}$ and
establish the following
\begin{prop}
The  following expressions for  the minors of $x w_0^{-1}$ in terms of the
variables $y_{k,i}$ and $T_{k,i}$ hold:
\be\label{gg}
x_{i,i+1}=\sum_{k=i}^{N-1}
y_{k,i}=\sum_{k=i}^{N-1}e^{T_{k,i}-T_{k+1,i+1}+T_{N,i+1}-T_{N,i}},\\
\Delta_i(x w_0^{-1})=(-1)^{i-1} \prod_{k=i}^{N-1}\prod_{j=k-i+1}^ky_{k,j}=
(-1)^{i-1} e^{\sum_{k=1}^i T_{i,k}-T_{N,k}},\\
\frac{\Delta_{i,i+1}(x w_0^{-1})}{\Delta_i(xw_0^{-1})}=
\sum_{k=i}^{N-1}\prod_{m=k}^{N-1}\frac{y_{m+1,m-i+1}}{y_{m,m-i+1}}=
\sum_{k=1}^{N-i}e^{T_{k+i,k}-T_{k+i-1,k}},
\ee
where we assume  $y_{N,k}=1$.
\end{prop}
Using the above relations it is easy to see that the integral
representation given in $\,\,\,$ \cite{GKMMO} reduces to the Gauss-Givental
representation. The direct between the two integral representations
gives also an explicit definition to the integration contour used in \cite{GKMMO}.

\section{ The Baxter $Q-$operator induction over  the rank}

\subsection{Induction over the rank and the Gauss-Givental construction}

The easiest way to construct the  Whittaker vectors in the
Gauss-Givental parameterization uses the recursive structure of this
parameterization. It is shown in Appendix B that using the conjugation 
by the  operators
\be \Xi_L^{(n)}=
e^{-\sum\limits_{i=1}^{n}\mu^{(n)}_i T_{n,i}}\exp\left\{
\mu^{(n)}_n\left(\sum_{i=1}^{n}T_{n,i}-\sum_{i=1}^{n-1}T_{n-1,i}\right)+\right.\\
\left.\sum_{i=1}^{n-1}\xi_L^{(n-i)}e^{T_{n,i}-T_{n-1,i}}\right\}
e^{\sum\limits_{i=1}^{n-1}\mu^{(n)}_i T_{n-1,i}}, \ee and
\be \Xi^{(n)}_R\,=
\exp\,\left\{-\sum_{i=1}^{n-1}T_{n,i}\sum_{k=1}^{n-1}
\frac{\partial}{\partial T_{k,i}}+
\sum_{i=2}^{n}T_{n,i}\sum_{k=1}^{i-1} \frac{\partial}{\partial
T_{k+(n-i),k}}\right\}\times\\
\exp\left\{\sum\limits_{i=1}^{n-1}\xi_R^{(i)}e^{T_{n-1,i}-T_{n,i+1}}\right\}
\times\\
\exp\left\{\sum_{i=1}^{n-2}T_{n-1,i}\sum_{k=1}^{n-2}
\frac{\partial}{\partial T_{k,i}}-
\sum_{i=2}^{n-1}T_{n-1,i}\sum_{k=1}^{i-1} \frac{\partial}{\partial
T_{k+(n-1-i),k}}\right\},\ee 
one can obtain the  iterative representations for the Whittaker
vectors in the form 
\be \psi_R^{(N)}=e^{H_L^{(N)}}\Xi^{(N)}_R \cdots \Xi^{(2)}_R\cdot 1, \ee
\be \psi_L^{(N)}=e^{H_L^{(N)}}\Xi^{(N)}_L\cdots \Xi^{(2)}_L\cdot 1. \ee
Let  us define \be  R= e^{-\sum_i T_{N,i}\rho_i^{(N)}}
\cdot\prod_{i=1}^{N-1}e^{(N-i)(T_{N-1,i}-T_{N,i})}.\ee
Then by  direct calculation  we can prove  the following recursive
 relation for  the Cartan elements 
\bqa
 R \cdot\overline{
\Xi^{(N)}_L}e^{-\sum T_{N,i}E^{(N)}_{i,i}}\Xi^{(N)}_R= e^{-\sum
T_{N-1,i}\rho^{(N-1)}_i}Q^{(N)}_{\lambda_N}\cdot e^{-\sum
T_{N-1,i}E^{(N-1)}_{i,i}}, \eqa 
where 
\be Q^{(N)}_{\lambda_N}(T_{N,1},\ldots,T_{N,N};T_{N-1,1},\ldots, T_{N-1,N-1})=\\
\hspace{-0.3cm}\exp\left\{\frac{\i\lambda_N}{\hbar}\Big(\sum_{k=1}^N T_{N,k}-
\sum_{k=1}^{N-1}T_{N-1,k}\Big)-\frac{1}{\hbar}
\sum_{k=1}^{N-1}\Big(e^{T_{N,k}-T_{N-1,k}}+e^{T_{N-1,k}-T_{N,k+1}}\Big)
\right\}.\ee
Taking into account these recursive relations   we obtain  
\begin{prop}\label{recrel}
\be
\Psi^{(N)}_{\la_1,\ldots,\la_N}(T_{N,1},\ldots,T_{N,N})=\\ \hspace{-2mm}
\int\limits_{\RR^{N-1}}
Q^{(N)}_{\lambda_N}(T_{N,1},\ldots,T_{N-1,N-1})
\Psi^{(N-1)}_{\la_1,\ldots,\la_{N-1}}(T_{N-1,1},\ldots,T_{N-1,N-1})
\prod\limits_{i=1}^{N-1}dT_{N-1,i}\,.
\ee
\end{prop}
{\it Proof.} Direct calculation.

Thus, the integral operator ${\cal Q}^{(N)}_{\lambda_N}$ transforms the
common eigenfunctions of $N-1$ particle open Toda chain Hamiltonians
to the ones for $N$ particle open Toda chain Hamiltonians.
Therefore, one can rederive Givental's integral formula by applying iteratively integral
operators ${\cal Q}^{(n)}_{\lambda_n}$, $1<n\leq N$. In the next subsection
we provide a more conceptual proof of the Proposition \ref{recrel}
and discuss the connection with Baxter's $Q$-operator formalism.

\subsection{Recursive derivation and  Baxter $Q-$operator}
Surprisingly, the Gauss-Givental representation has a direct connection
with the formalism of the recursion operator \cite{KL1} and thus a close
relation with the Baxter $Q$-operator and the QISM approach. Below we rederive
this integral representation using the recursion procedure. Let us first
recall the formalism of the Lax operators (see e.g. \cite{Ga}). In this
subsection we set $\xi^{(i)}_R=\xi^{(i)}_L=-\h^{-1}$ for simplicity.

Consider the following $N\times N$ Lax operators
\be\label{tc1}
L_N=\begin{pmatrix}
p_1&1&0&0&...&0\\
 e^{x_1-x_2}&p_2&1&0&...&0\\
0&...&...&...&...&\\
...&...&...&...&...&\\
0&0&...& e^{x_{N-2}-x_{N-1}}&p_{N-1}&1\\
0&0&...&0& e^{x_{N-1}-x_N}&p_N
\end{pmatrix},
\ee
where $p_{n}=-\i\hbar \frac{\pr}{\pr x_{n}}$. The generating function for the
Hamiltonians $H_n^{(N)}$ of the $N$ particle open Toda chain is defined as
\be\label{tc3}
\det(u-L_N)=\sum_{n=0}^N(-1)^nu^{N-n}H_n^{(N)}\,.
\ee
The determinant $A_N(u)\equiv \det(u-L_N)$ with non-commutative elements is
defined iteratively
\be
A_n(u)=(u-p_n)A_{n-1}(u)-e^{x_{n-1}-x_n}A_{n-2}(u),
\ee
for $n=1,2,\ldots$, where  $A_{-1}(u)=0$, $A_0(u)=1$.

For any function $f\in{\rm Fun}\,\RR^{N-1}$ not growing too fast at infinity,
let
\be\label{in1}
({\cal Q}^{(N)}_u f)(x_1,\ldots,x_N)=\\=\int\limits_{\RR^{N-1}}
Q^{(N)}_u (x_1,\ldots,x_N;y_1,\ldots,y_{N-1})f(y_1,\ldots,y_{N-1})
dy_1\ldots dy_{N-1}\,,
\ee
where
\be\label{in2}
Q^{(N)}_u (x_1,\ldots,x_N;y_1,\ldots,y_{N-1}):=\\=
\exp\left\{\frac{\i u}{\h}\Big(\sum_{i=1}^Nx_i-\sum_{i=1}^{N-1}y_i\Big)-
\frac{1}{\h}\sum_{i=1}^{N-1}\Big(e^{x_i-y_i}+e^{y_i-x_{i+1}}\Big)\right\}.
\ee
\begin{lem}
The following intertwining relation holds:
\be\label{in3}
A_N(u)\circ {\cal Q}^{(N)}_v=(u-v){\cal Q}^{(N)}_v \circ A_{N-1}(u)\,,
\ee
where $A_{N-1}(u):=\sum\limits_{n=0}^{N-1}(-1)^nu^{N-n-1}H^{(N-1)}_n$
is the generating function for the Hamiltonians of the $N-1$ particle open
Toda chain acting as a differential operators in variables $y_1,\ldots,y_{N-1}$.
\end{lem}
{\it Proof.} The proof is similar to the proof of Theorem 2 in \cite{Gi}
(see also  \cite{JK}).

\square
\begin{prop}
Let $\Psi_{\la_1,\ldots,\la_{N-1}}^{(N-1)}(T_{N-1,1},\ldots,T_{N-1,N-1})$ be
the eigenfunction  of the $N-1$ particle open Toda chain. Then the function
\be\label{in6}
\Psi_{\la_1,\ldots,\la_N}^{(N)}(T_{N,1},\ldots,T_{N,N})=
({\cal Q}^{(N)}_{\la_N}\Psi_{\la_1,\ldots,\la_{N-1}}^{(N-1)})(T_{N,1},
\ldots,T_{N,N})
\ee
is the solution of the $N$ particle open Toda chain depending on  coordinates
$T_{N,1},\ldots,T_{N,N}$.\\
Here the integral kernel of the  operator is defined by
\be\label{k2}
Q^{(N)}_{\lambda_N} (T_{N,1},\ldots,T_{N,N};
T_{N-1,1}\ldots,T_{N-1,N-1})=\\ \hspace{-0.6cm}
\exp\left\{\frac{\i\lambda_N}{\hbar}\Big(\sum_{i=1}^N T_{N,i}-
\sum_{i=1}^{N-1}T_{N-1,i}\Big)-\frac{1}{\hbar}
\sum_{i=1}^{N-1}\Big(e^{T_{N,i}-T_{N-1,i}}+e^{T_{N-1,i}-T_{N,i+1}}\Big)
\right\}.
\ee
\end{prop}
{\it Proof}. Assume that the function
$\Psi_{\la_1,\ldots,\la_{N-1}}^{(N-1)}(T_{N-1,1},\ldots,T_{N-1,N-1})$
satisfies the equation
$
A_{N-1}(u)\Psi_{\la_1,\ldots,\la_{N-1}}^{(N-1)}=
\prod_{k=1}^{N-1}(u-\la_k)\Psi_{\la_1,\ldots,\la_{N-1}}^{(N-1)}
$.
Using (\ref{in3}) with $u=\la_N$, one obtains that the function
(\ref{in6}) satisfies the equation $A_N(u)\Psi_{\la_1,\ldots,\la_N}^{(N)}=
\prod_{k=1}^N(u-\la_k)\Psi_{\la_1,\ldots,\la_N}^{(N)}$.
\square

\begin{cor}
The integral representation (\ref{Givp}), (\ref{pot}) can be obtained by
iteration of the recursive relation (\ref{in6}).
\end{cor}

Finally note that the recursive operator $\mathcal{Q}$ introduced above is
closely connected with the Baxter $Q$-operator of the periodic Toda chain
\cite{PG}. The kernel of the operator  $\mathcal{Q}$ is obtained form the
kernel of the Baxter $Q$-operator by dropping out the terms corresponding to
the ``interaction'' of the first and last nodes of the Toda chain.
Moreover, the elementary blocks of the parameterization of the unipotent
matrices introduce in Section 2 also have counterparts in the
Baxter $Q$-operator formalism (see e.g. eq.(54) in  \cite{Sk2}).

\section{On a generalization to infinite and periodic Toda chain}

The Gauss-Givental representation of the open Toda chain wave function
seems to imply  a  natural generalization to the case of the infinite
and periodic Toda chain wave function. The details of this
generalizations will be discussed elsewhere.
 In this section we discuss the formal expressions arising in this case. For simplicity we
will consider below only the case of zero eigenvalues  $\lambda_i=0$ of the
Hamiltonians.

  Let us start with the infinite  Toda chain. The quantum
Hamiltonians of the infinite  Toda chain are obtained by taking the formal
 limit $N\rightarrow \infty$  of the finite open Toda case  (see for
 example \cite{Sk1}).
  Thus the first nontrivial Hamiltonian is given by
\be
H=-\frac{\hbar^2}{2}\sum\limits_{i=-\infty}^{\infty}\frac{\partial^2}{\partial x_i^2}+
\sum\limits_{i=-\infty}^{\infty}e^{x_{i}-x_{i+1}}.
\ee
It is natural to guess that the generalization of the integral
representation for the common eigenfunctions of the commuting set of
 Hamiltonians for the infinite Toda chain is given by the appropriately
defined limit $N\rightarrow \infty$ of the Gauss-Givental integral
representation
\be\label{Givinf}
\Psi^{(\infty)}(T_{0,i})=
\int_{\Gamma}
e^{\frac{1}{\hbar}\mathcal{F}^{(\infty)}(T)}\prod_{k=1}^\infty
\prod_{i=-\infty}^\infty dT_{k,i}\,,
\ee
where the function $\mathcal{F}^{(\infty)}(T)$ is given by
\bqa\label{potinf}
\mathcal{F}^{(\infty)}(T)=-\sum_{k=1}^{\infty}
\left(\sum_{i=-\infty}^{\infty}e^{T_{k-1,i}-T_{k,i}}+
\sum_{i=-\infty}^{\infty}e^{T_{k,i}-T_{k-1,i+1}}\right),
\eqa
and the cycle  $\Gamma$  is a semiinfinite-dimensional   submanifold in  the
infinite-dimensional  complex torus with coordinates
$\{e^{T_{k,j}}|\,\,\, k\in \mathbb{N},i\in \mathbb{Z}\}$. Note that we
 relabel  the variables $T_{k,i}$ to take
the limit $N\rightarrow \infty$ properly.

Similarly one defines the formal integral representation for the
periodic Toda chain. It is well known that the Hamiltonians of the
$N$ particle periodic Toda chain  can
be obtained from the Hamiltonians of the infinite Toda chain by
imposing the periodicity condition $x_{i+N}=x_i$. Thus the first
non-trivial Hamiltonian is given by
\be
H=-\frac{\hbar^2}{2}\sum\limits_{i=1}^{N}\frac{\partial^2}{\partial x_i^2}+
\sum\limits_{i=1}^{N}e^{x_{i}-x_{i+1}},
\ee
and $x_{N+1}=x_1$. It is natural to guess that similar factorization works
on the level of the integral representation. Thus  we get formally
\be\label{Givper}
\Psi^{(per)}(T_{0,i})=\int_{\Gamma}
e^{\frac{1}{\hbar}\mathcal{F}^{(per)}(T)}\prod_{k=1}^{\infty}
\prod_{i=1}^N dT_{k,i}\,\,,
\ee
where the function $\mathcal{F}_{(per)}(T)$ is given by
\bqa\label{potper}
\mathcal{F}^{(per)}(T)=-\sum_{k=1}^{\infty}
\left(\sum_{i=1}^{N}e^{T_{k-1,i}-T_{k,i}}+
\sum_{i=1}^{N}e^{T_{k,i}-T_{k-1,i+1}}\right),
\eqa
where we assume $T_{k,N+i}=T_{k,i}$ and the cycle  $\Gamma$ is a semi-infinite
dimensional submanifold in the infinite-dimensional  complex torus with
coordinates $\{e^{T_{k,j}}|\,\,\, k\in \mathbb{N},i=1,\ldots ,N\}$. Although
these integral representations are rather meaningless, they seem reveal the
important properties of the solution. For example the infinite number of
integrations in the periodic case is a manifestation of the integration over
the ``positive part'' of the loop group $LGL(N,\mathbb{R})$ which is a formal
analog of the unipotent upper triangular matrices in the affine case.

Let us note that the formal expressions (\ref{Givinf}) and (\ref{Givper}) have
an interesting interpretation from the point of view of the recursive
relations discussed in the previous section. Consider the formal limit
$N\rightarrow \infty$ of the integral kernel of the recursion operator in the
case of the infinite Toda chain
\be\label{kinf}
Q^{(\infty)}(T_{k,*},T_{k-1,*})=\exp\left\{-\frac{1}{\hbar}
\sum_{i=-\infty}^{\infty}\Big(e^{T_{k-1,i}-T_{k,i}}+e^{T_{k,i}-T_{k-1,i+1}}
\Big)\right\}.
\ee
It is easy to see that the integral operator defined by this kernel function
commutes with the Hamiltonians and  is a particular example of the
B\"{a}cklund transformation. This kind of transformations for infinite Toda
chain defined by the integral kernel (\ref{kinf}) was studied previously
in detail \cite{T}. Similarly in the case of the periodic Toda chain we have
that
\be\label{kinper}
Q^{(per)}(T_{k,1},\ldots,\,T_{k,N},T_{k-1,1}\ldots,
\,T_{k-1,N}) =\\ \exp\left\{-\frac{1}{\hbar}
\sum_{i=1}^{N}\Big(e^{T_{k-1,i}-T_{k,i}}+e^{T_{k,i}-T_{k-1,i+1}}\Big)
\right\},
\ee
where we  assume $T_{k,i+N}=T_{k,i}$ and $T_{k-1,i+N}=T_{k-1,i}$.
As in the case of the infinite chain the corresponding integral
transformation 
provides an example of a  B\"{a}cklund transformation
for the periodic chain. Exactly these transformations were considered
previously in \cite{Ga}, \cite{PG}, \cite{Sk2}. Taking into account the
explicit expressions (\ref{kinf}), (\ref{kinper}) for the kernels the integral
formulas (\ref{Givinf}) and (\ref{Givper}) can  be also interpreted as
follows. Let us apply $K$ times the integral operator with the kernel
(\ref{kinf}) to the constant function. It is easy to see that when
$K\rightarrow \infty$ the resulting expression tends to the  one given by the
integral formula (\ref{Givinf}). This is consistent with the fact that the
action of $Q^{(\infty)}$ transforms an  eigenfunction of the infinite Toda
chain Hamiltonians into another eigenfunction. Completely similar picture
holds for the periodic case.

Finally let us mention that the solution for the common eigenfunction of the
Hamiltonians of the periodic Toda chain has already appeared elsewhere
\cite{KL2}, \cite{KL3}. It is given by the finite-dimensional integral with
rather complex integrand. Hopefully it could be obtained after an infinite
number of the integrations from an appropriately defined formal integral
representation (\ref{Givper}).

\section{Conclusions}
In this paper we have proposed an interpretation  of the solution of the open
Toda chain given in \cite{Gi} (see also \cite{JK}) in terms of representation
theory. This leads to a particular realization of the principal series
representation of $U(\mathfrak{gl}(N))$ in terms of the differential operators
acting on the space of functions defined on the subspace of the totally
positive unipotent upper-triangular matrices.  The resulting representation
is very close to the representation constructed (rather implicitly) in
\cite{BFZ} which is based on the standard factorization of the unipotent
upper-triangular matrices into the product of the elementary Jacobi matrices.
Note that the representation of \cite{BFZ} can be generalized to the case of
an arbitrary semisimple Lie algebra. Explicitly the parameterization is given
in terms of the matrix elements of the fundamental representations
(see e.g \cite{Li}). This implies that the Gauss-Givental representation
discussed in this paper may be also rather straightforwardly generalized to
the case of an arbitrary semisimple Lie algebra. Let us remark that the
concept of the total positivity was generalized to the case of an arbitrary
semisimple group by  Lusztig (see \cite{Lu} for details) in connection with
the study of canonical bases. It is interesting that the same concept seems
to play an important role in the choice of the proper coordinates for the
construction of the Gauss-Givental representation proposed in this paper.

The connection between the Baxter $Q$-operator and the recursion operator
discussed in Section 4.2 provides a simple  construction of the Givental type
integral representations for the common eigenfunctions of Hamiltonians for a
wide class of the models (e.g. XXX and XXY open chains). These applications
will be discussed elsewhere \cite{GKLO2}. Also note that the integral
representations of the wave function of the periodic Toda chain corresponding
to the affine algebra $\hat{\mathfrak{gl}}_N$ provides an interesting
description of the correlators in the Landau-Ginzburg topological model.
Mirror symmetry predicts that on the other hand one gets the description of
the Gromov-Witten invariants of the affine flag manifold.

Finally one should say that there is a generalization of the results of this
paper to the case of the quantum groups and  the centrally extended loop
algebras. We are going to discuss various generalization of the Gauss-Givental
representation constructed in this paper elsewhere.

\appendix

\section{Appendix A}

In this appendix we give the proof of the Proposition \ref{baserep}.
 We will  organize our  calculations in two steps. First we
use  the recursive  structure of (\ref{factor})  and obtain  recursive formulas for
 generators  $E^{(N)}_{i+1,i}$, $E^{(N)}_{i,i}$ and $E^{(N)}_{i,i+1}$
 in terms of $y_{k,i}$. Then,  changing  variables and rewriting  the recursive relations
in variables $T_{k,i}$,  we solve them to  obtain the final form of
the generators.

\subsection{The case of $E^{(N)}_{i,i+1}$}
 In order to find an  explicit form of the  generators $E_{i,i+1}^{(N)}$
in the variables $y_{i,k}$ we use the following commutation
relations between elementary unipotent matrices:
\be\label{r1}
\left(1+\sum_{i=1}^{N-1}y_{N-1,i}\cdot
e_{i,i+1}\right)\left(1+\varepsilon\cdot e_{j,j+1}\right)=\\
\left(1+\frac{y_{N-1,j-1}}{y_{N-1,j}}\varepsilon\cdot e_{j-1,j}\right)
\left(1+\sum_{i=1}^{N-1}y'_{N-1,i}\cdot e_{i,i+1}\right) mod(\varepsilon^2)
,\ee
where
\be
y'_{N-1,j-1}=y_{N-1,j-1}-\frac{y_{N-1,j-1}}{y_{N-1,j}}\varepsilon ,\\
\,\, y'_{N-1,j}=y_{N-1,j}+\varepsilon,\,\,
y'_{N-1,i}=y_{N-1,i},\,i\neq j-1,j. \ee
Using  the  decomposition  (\ref{factor}),
definition (\ref{def}) and the relation (\ref{r1})
it is  simple to check that, in variables $y_{i,j}$, the  generators
$E^{(N)}_{i,i+1}$ obey the following recursion relation:
\be
E^{(N)}_{i,i+1}f(x)=\\ \left(\frac{y_{N-1,i-1}}{y_{N-1,i}}E^{(N-1)}_{i-1,i}-
y_{N-1,i}^{-1}\left(y_{N-1,i-1}\frac{\partial}{\partial
y_{N-1,i-1}}+y_{N-1,i}\frac{\partial}{\partial y_{N-1,i}}\right)\right)f(x).
\ee
Resolving the recursion  procedure we obtain:
\be
\hspace{-0.6cm}
E^{(N)}_{i,i+1}=\sum_{k=0}^{i-1}\prod_{s=0}^k
\frac{y_{N-s,i-s}}{y_{N-s,i+1-s}}\frac{\partial}{\partial
y_{N-1-k,i-k}}- \prod_{s=0}^k\frac{y_{N-(s+1),i-(s+1)}}{y_{N-(s+1),i-s}}
\frac{\partial}{\partial y_{N-1-k,i-(k+1)}}=\\
\sum_{k=1}^i\left(\,\sum_{s=k}^ie^{T_{N+s-i,s}-T_{N+s-i-1,s}}\right)
\times\left(\frac{\partial}{\partial T_{N+k-i-1,k}}-
\frac{\partial}{\partial T_{N+k-i-1,k-1}}\right).
\ee
Here we have used:
\bqa\label{ch2}
y_{p,n}\frac{\partial}{\partial y_{p,n}}=
\sum\limits_{k=1}^{n}\frac{\partial}{\partial T_{p-k+1,n-k+1}},\,\,\,\,\,
n=1,\ldots,p.
\eqa

\subsection{The case of $E^{(N)}_{i,i}$}

For the Cartan generators we get the same type of recursion
procedure in the variables $y_{k,i}$: \bqa
E^{(N)}_{i,i}f(x)&=&\left(E^{(N-1)}_{i,i}+y_{N-1,i-1}\frac{\partial}{\partial
y_{N-1,i-1}}-y_{N-1,i}\frac{\partial}{\partial
y_{N-1,i}}\right)f(x),\,\,\,\ i\neq N,\\
E^{(N)}_{N,N}f(x)&=&(\mu_N^{(N)}+y_{N-1,N-1}\frac{\pr}{\pr y_{N-1,N-1}})f(x),
\eqa
using  the following relation:
\be
\hspace{-0.5cm}
\left(1+\sum_{i=1}^{N-1}y_{N-1,i}\cdot
e_{i,i+1}\right)\left(1+\varepsilon\cdot e_{j,j}\right)=
(1+\varepsilon\cdot e_{j,j})
\left(1+\sum_{i=1}^{N-1}y'_{N-1,i}\cdot e_{i,i+1}\right),
\ee
where \be y'_{N-1,j-1}=y_{N-1,j-1}(1+\varepsilon) , \qquad
y'_{N-1,j}=y_{N-1,j}(1+\varepsilon)^{-1},\qquad\\
\widetilde{y}_{N-1,i}=y_{N-1,i},\,i\neq j,j-1. \ee
Resolving recursion relations one can find
\bqa \label{cart}
E^{(N)}_{i,i}=
\mu^{(N)}_i + \sum_{k=1}^{i-1}\frac{\partial}{\partial
T_{N+k-i,k}} - \sum_{k=i}^{N-1}\frac{\partial}{\partial
T_{k,i}}.
\eqa

\subsection{The case of $E^{(N)}_{i+1,i}$}

In this case the same line of reasoning leads to  the following
recursion
expression:
\be E^{(N)}_{i+1,i}f(x)=
\left(E^{(N-1)}_{i+1,i}+y_{N-1,i}\left(
E^{(N-1)}_{i,i}-E^{(N-1)}_{i+1,i+1}\right)\right.-\\
y_{N-1,i}\left.\left(y_{N-1,i}\frac{\partial}{\partial
y_{N-1,i}}+y_{N-1,i+1}\frac{\partial}{\partial
y_{N-1,i+1}}\right)\right) f(x).
\ee
Resolving the recursion procedure and using (\ref{cart}) one obtains:
\be
E^{(N)}_{i+1,i}=\sum_{k=1}^{N-1}\left[(\mu^{(N)}_i-\mu^{(N)}_{i+1})
y_{k,i}-y_{k,i}(y_{k,i}\frac{\partial}{\partial
y_{k,i}}-y_{k,i+1}\frac{\partial}{\partial y_{k,i+1}})+\right.\\
\left.y_{k,i}\sum_{s=1}^{k-1}(y_{s,i-1}\frac{\partial}{\partial
y_{s,i-1}}-2y_{s,i}\frac{\partial}{\partial y_{s,i}}+
y_{s,i+1}\frac{\partial}{\partial y_{s,i+1}})\right]=\\
\sum_{k=1}^{N-1}e^{(T_{k,i}-T_{k+1,i+1})}
\left(\mu^{(N)}_i-\mu^{(N)}_{i+1}+\sum_{s=1}^k\left(\frac{\partial}{\partial
T_{s,i+1}}-\frac{\partial}{\partial T_{s,i}}\right)\right).
\ee
\square

Notice that the generators $E^{(N)}_{i,i}, E^{(N)}_{i,i+1}, E^{(N)}_{i+1,i}$ written
as differential operators in terms
of $y_{k,i}$  has a natural recursive
structure.  However  the variables $T_{k,i}$ are more appropriate
when the   problem of solving  the equations on Whittaker vectors
arises.

\section{Appendix B}

In this Appendix we give  the proof of the Proposition \ref{Wvec}.
To derive the explicit  expressions for the Whittaker vectors
$\psi_{L,R}$  we use the iteration over the rank of the algebra
$\mathfrak{gl}(N)$. Consider the twisted generators of
the algebra $\mathfrak{gl}(N)$ 
\bqa\label{cartnew} \widehat{E}^{(N)}_{i,i}&=&\mu^{(N)}_i+
\sum_{k=1}^{i-1}\frac{\partial}{\partial
T_{N+k-i,k}}-\sum_{k=i}^{N-1}\frac{\partial}{\partial
T_{k,i}},\\
\label{rasnew}
\widehat{E}^{(N)}_{i,i+1}&=&\sum_{n=1}^i\left(\,\sum_{k=n}^ie^{T_{N+k-i,k}-T_{N+k-i-1,k}}\right)
\left(\frac{\partial}{\partial T_{n+N-i-1,n}}-
\frac{\partial}{\partial T_{n+N-i-1,n-1}}\right),\\
\widehat{E}^{(N)}_{i+1,i}&=&\sum_{k=1}^{N-1}e^{(T_{k,i}-T_{N,i})-(T_{k+1,i+1}-T_{N,i+1})}
\left(\mu^{(N)}_i-\mu^{(N)}_{i+1}+\sum_{s=1}^k\left(\frac{\partial}{\partial
T_{s,i+1}}-\frac{\partial}{\partial T_{s,i}}\right)\right), \eqa 
that are related with the generators (\ref{baserep}) as follows
\be 
\widehat{E}_{i,j}^{(N)}=e^{-H_L^{(N)}}E_{i,j}^{(N)}e^{H_L^{(N)}}. \ee
To reveal the recursive structure it is useful to   modify the
generators $\widehat{E}^{(N)}_{i,i+1}$ and
$\widehat{E}^{(N)}_{i+1,i}$ separately 
 \bqa\label{Dy}
\widetilde{E}^{(N)}_{i,i+1}:=(\Xi^{(N)}_R)^{-1}\widehat{E}^{(N)}_{i,i+1}
(\Xi^{(N)}_R)\eqa \bqa\widetilde{
\widetilde{E}}^{(N)}_{i+1,i}:=(\Xi^{(N)}_L)^{-1}\widehat{E}^{(N)}_{i+1,i}
(\Xi^{(N)}_L)\eqa where the operators $\Xi^{(N)}_{L,R}$ are defined
as   \be\label{XiL} \Xi_L^{(n)}=
e^{-\sum\limits_{i=1}^{n}\mu^{(n)}_i T_{n,i}}\exp\left\{
\mu^{(n)}_n\left(\sum_{i=1}^{n}T_{n,i}-\sum_{i=1}^{n-1}T_{n-1,i}\right)+\right.\\
\left.\sum_{i=1}^{n-1}\xi_L^{(n-i)}e^{T_{n,i}-T_{n-1,i}}\right\}
e^{\sum\limits_{i=1}^{n-1}\mu^{(n)}_i T_{n-1,i}}, \ee and
\be\label{XiR} \Xi^{(n)}_R\,=
\exp\,\left\{-\sum_{i=1}^{n-1}T_{n,i}\sum_{k=1}^{n-1}
\frac{\partial}{\partial T_{k,i}}+
\sum_{i=2}^{n}T_{n,i}\sum_{k=1}^{i-1} \frac{\partial}{\partial
T_{k+(n-i),k}}\right\}\times\\
\exp\left\{\sum\limits_{i=1}^{n-1}\xi_R^{(i)}e^{T_{n-1,i}-T_{n,i+1}}\right\}
\times\\
\exp\left\{\sum_{i=1}^{n-2}T_{n-1,i}\sum_{k=1}^{n-2}
\frac{\partial}{\partial T_{k,i}}-
\sum_{i=2}^{n-1}T_{n-1,i}\sum_{k=1}^{i-1} \frac{\partial}{\partial
T_{k+(n-1-i),k}}\right\}.\ee Here $n=1,\ldots,N$ and we also set
$\Xi^{(1)}_{L,R} =1$. 
We have the following relations 
\be\label{tld}
e^{T_{N-1,i}-T_{N,i+1}}\cdot\left(\,\widetilde{E}^{(N)}_{i,i+1}-\xi^{(i)}_R\,\right)\,=\\
e^{T_{N-1,i-1}-T_{N,i}}\cdot
\left(\,\widehat{E}^{(N-1)}_{i-1,i}-\xi^{(i-1)}_R\,\right)\,+\,
\left(\,\sum_{k=i}^{N-1}\frac{\partial}{\partial
T_{k,i}}-\sum_{k=i-1}^{N-1}\frac{\partial}{\partial
T_{k,i-1}}\,\right)\ee
and \be\label{dtld}
e^{T_{N,i}-T_{N-1,i}}\cdot\left(\,\widetilde{\widetilde{E}}^{(N)}_{i+1,i}-\xi^{(i)}_L\,\right)\,=\\
e^{T_{N,i+1}-T_{N-1,i+1}}\cdot
\left(\,\widehat{E}^{(N-1)}_{i+1,i}-\xi^{(i-1)}_L\,\right)\,-\,
\left(\,\sum_{k=i}^{N-1}\frac{\partial}{\partial
T_{k,i}}-\sum_{k=i+1}^{N-1}\frac{\partial}{\partial
T_{k,i+1}}\,\right).\ee
\begin{lem}
The  equations 
\be
\hspace{2cm}
\widehat{E}^{(N)}_{i,i+1}\widehat{\psi}^{(N)}_R=\xi^{(i)}_R \widehat{\psi}^{(N)}_R\,,
\hspace{1cm}(i=1,\ldots,N-1),
\ee
\be
\hspace{2cm} \widehat{E}^{(N)}_{i+1,i} \widehat{\psi}^{(N)}_L=\xi^{(N-i)}_L \widehat{\psi}^{(N)}_L\,,
\hspace{1cm} (i=1,\ldots,N-1).
\ee
 admit the following solution   \be\label{itR} \widehat{\psi}_R^{(N)}=\Xi^{(N)}_R
\widehat{\psi}_R^{(N-1)}=\Xi^{(N)}_R \cdots \Xi^{(2)}_R\cdot 1, \ee
\be\label{itL} \widehat{\psi}_L^{(N)}=\Xi^{(N)}_L
\widehat{\psi}_L^{(N-1)}=\Xi^{(N)}_L\cdots \Xi^{(2)}_L\cdot 1. \ee
\end{lem}
{\it Proof}. The statement follows  from the recursion
representation (\ref{tld}) and (\ref{dtld}) for the generators.
\square
Now the Whittaker vectors (\ref{wR}) and (\ref{wL})
are obtained as follows
\be
\widehat{\psi}^{(N)}_{L,R}=e^{-H^{(N)}_L}\psi^{(N)}_{L,R}.\ee

\end{document}